\documentclass[10pt,a4paper]{amsart}
\usepackage{amssymb, amsthm, amsmath}

\newtheorem{thm}{Theorem}[section]
\newtheorem{lem}[thm]{Lemma}
\newtheorem{prop}[thm]{Proposition}
\newtheorem{cor}[thm]{Corollary}
\theoremstyle{definition}
\newtheorem{defn}[thm]{Definition}

\newcommand{\h}{\mathcal{H}}
\newcommand{\dir}[1]{\oplus_{i\in\Omega}#1}
\newcommand{\dirk}[1]{\oplus^k#1}
\newcommand{\normi}[1]{\|#1\|_i^2}
\newcommand{\norm}[1]{\|#1\|^2}
\newcommand{\iri}{E_\lambda^i}
\newcommand{\ltwoi}{L^2(M_i,\mu_i)}
\newcommand{\normiri}[1]{\|E_\lambda^i#1\|^2}
\newcommand{\intinf}[2]{\int_{-\infty}^\infty #1\,d#2}

\begin{document}

\title{Introduction to the spectral theory of self-adjoint differential vector-operators}%
\author{Maksim S. Sokolov}
\address{Maksim S. Sokolov \\ ICTP Affiliated Center, Mechanics and Mathematics Department, National University of Uzbekistan (Uzbekistan, Tashkent 700095)}%
\email{sokolovmaksim@hotbox.ru}%

\maketitle

\begin{abstract}
We study the spectral theory of operators, generated as direct
sums of self-adjoint extensions of quasi-differential minimal
operators on a multi-interval set (self-adjoint vector-operators),
acting in a Hilbert space.

Spectral theorems for such operators are discussed, the structure
of the ordered spectral representation is investigated for the
case of differential coordinate operators. One of the main results
is the construction of spectral resolutions. Finally, we study the
matters connected with analytical decompositions of generalized
eigenfunctions of such vector-operators and build a matrix
spectral measure leading to the matrix Hilbert space theory.
Results, connected with other spectral properties of
vector-operators, such as the introduction of the identity
resolution and the spectral multiplicity have also been obtained.

Vector-operators have been mainly studied by W.N. Everitt, L.
Markus and A. Zettl. Being a natural continuation of
Everitt-Markus-Zettl theory, the presented results reveal the
internal structure of self-adjoint vector-operators and are
essential for the further study of their spectral properties.
\end{abstract}

\section{Preliminaries}
\subsection{Problem overview}
In 1985, F. Gesztesy and W. Kirsch published their work
\cite{shrodinger}, where they considered an example of a
Schr\"odinger operator generated by the Hamiltonian of the type
\begin{equation}\label{hamiltonian}
  H=-\frac{d^2}{dx^2}+ \frac1{\cos^2x}.
\end{equation}
Since the potential of (\ref{hamiltonian}) has a countable number
of singularities on $\mathbb{R}$ which spoil the local
integrability, they constructed minimal operators $T_{i,min}$,
generated by (\ref{hamiltonian}) in the spaces
$$L^2\left(-\frac\pi 2+i\pi,\frac\pi
2+i\pi\right),\,\,i\in\mathbb{Z},$$ and then considered the direct
sum operator $\oplus_{i\in\mathbb{Z}}T_{i,min}$ in the space
\begin{equation*}\label{space for sh op}
    \oplus_{i\in\mathbb{Z}}L^2\left(-\frac\pi 2+i\pi,\frac\pi 2+i\pi\right)
\end{equation*} which equaled to the minimal operator in
$L^2(-\infty,\infty).$

The work \cite{shrodinger} stimulated other researchers to
generalize the problem. In 1992, W.N.~Everitt and A.~Zettl
\cite{everittzettl} studied direct sums of minimal and maximal
operators generated by arbitrary formally self-adjoint expressions
in Hilbert spaces considered on arbitrary intervals (maximal and
minimal vector-operators). Later in 2000, vector-operators were
also considered in complete locally convex spaces by R.R.~Ashurov
and W.N.~Everitt \cite{ashuroveveritt}, which was a natural
generalization of their work \cite{ashuroveveritt1}. Since 1992,
quasi-differential vector-operators have mostly been investigated
in connection with their non-spectral theory, such as the
introduction of minimal and maximal vector-operators and their
relationship (it was shown that the adjoint of a minimal
vector-operator is maximal in a Hilbert space \cite{everittzettl},
and the analogous result with the modification for Frechet spaces
was obtained in \cite{ashuroveveritt}). A lot of work has been
carried out by W.N.~Everitt and L.~Markus in order to develop the
theory of self-adjoint extensions for vector-operators with the
employment of symplectic geometry. In connection with this, see
their recent memoirs \cite{everittmarkus} and
\cite{everittmarkus1}.

The theory of operators generated by multi-interval systems finds
its applications in many problems of quantum mechanics, theory of
semiconductors and theoretical computer science; good
bibliographical references for these subjects may be found in
\cite{everittmarkus1}.

Since the theory of quasi-differential vector-operators in a
Hilbert space is quite young and the most recent studies have
concerned mostly problems connected with their common theory,
small attention was given to its spectral aspects. Some results,
describing the position of spectra of vector-operators were
presented in 1985 in \cite{shrodinger} and the most recent results
belong to Sobhy El-Sayed Ibrahim \cite{ibrahim,ibrahim1}.

This work is devoted to the study of spectral properties of
self-adjoint vector-operators. The latter do not cover all the
possible self-adjoint extensions of minimal vector-operators and
they are physically less interesting, but, nevertheless, their
spectral theory seems to be mathematically interesting.

The abstract spectral theory of self-adjoint vector-operators was
presented in \cite{sokolov} and \cite{sokolov1}. Differential
coordinate operators play the key role in \cite{sokolov2} and
\cite{sokolov3}.

\subsection{Quasi-differential operators and vector-operators}
Basic concepts of quasi-differential operators are described in
\cite{everittzettl, everittmarkus}. A good reference for operators
with real coefficients is the book of M.A. Naimark \cite{naimark}.

Let us have a number $n\in \mathbb{N}, n\geqslant 2$, and an
arbitrary interval  $I\subseteq \mathbb{R}$. Let $Z_n(I)$ be a set
of Shin-Zettl matrices. These are matrices $A=\{a_{rs}\},\;
a_{rs}:I\to \mathbb{C}$ of the order $n\times n$, such that for
almost all $x\in I$:
\begin{equation*}\label{shinzetl}
  \left\{
  \begin{array}{lll}
  (i) & a_{rs}\in L_{loc}(I), & r,s = \overline{1,n}; \\{}\\ (ii) & a_{r,r+1}(x)\neq 0, & r= \overline{1,n-1};\\ {}\\(iii)&
   a_{rs} = 0, & s=\overline{r+2,n};\; r= \overline{1,n-2}.\\
  \end{array} \right.
\end{equation*}

Consider a function $f:I\to\mathbb{C}$; its quasi-derivatives
relative to a Shin-Zettl matrix $A$ are defined by
\begin{equation*}\label{quasideriv}
  \left\{
  \begin{array}{lll}
  (i) & f_A^{[0]}:=f; & {} \\{}\\ (ii) & f_A^{[r]}:=
  \frac 1{a_{r,r+1}}\left[ \frac d{dx}f_A^{[r-1]} - \sum_{s=1}^ra_{rs}f_A^{[s-1]}
  \right], & r= \overline{1,n-1};\\{}\\ (iii)&
f_A^{[n]}:= \frac d{dx}f_A^{[n-1]} -
\sum_{s=1}^na_{ns}f_A^{[s-1]}.
   \end{array} \right.
\end{equation*}

Let us introduce a linear manifold $D(A)\subset AC_{loc}(I)$:
\begin{equation*}\label{domaingeneral}
D_A(I):= \{f:I\to \mathbb{C}|\; f_A^{[r-1]}\in AC_{loc}(I)\quad
(r= \overline{1,n})\}.
\end{equation*}
It is possible to see, that $f\in D_A(I)$ implies $f_A^{[n]}\in
L_{loc}(I)$, and it is possible to prove that $D_A(I)$ is dense in
$L_{loc}(I)$.

Relative to a matrix $A\in Z_n(I)$, we have the quasi-differential
expression $M_A[f]=i^nf_{A}^{[n]}$, $f\in D_A(I)$.

The matrix $A^+\in Z_n(I)$ designates the Lagrange adjoint matrix
to $A$ if $
  A^+:=-L_n^{-1}A^*L_n,
$ where $A^*$ is the adjoint matrix, and $L_n=\{l_{rs}\}$ is the
$(n\times n)$-matrix, defined as:
$$l_{r,n+1-r}=\left\{
\begin{array}{cl}
(-1)^{r-1},&r= \overline{1,n};\\ 0, & \mbox{for other}\; r,s.
\end{array}\right.$$

Using this notation we suppose that in this work we deal only with
Lagrange symmetric (formally self-adjoint) expressions, that is
$M_{A^+}[f]=M_A[f]=\tau(f)$, where $\tau$ is an alternative
denotation for a Lagrange symmetric expression.

For a quasi-differential expression $M_A[f]$, the Lagrange formula
is known ($[\alpha,\beta]\subseteq I$ - an arbitrary compact
subinterval of $I$):
\begin{equation}\label{lagrangeformula}
\int_\alpha^\beta\{
\overline{g(x)}M_A[f](x)-f(x)\overline{M_{A^+}[g(x)]}\}\,dx=[f,g]_A(\beta)
- [f,g]_A(\alpha),
\end{equation}
where $f\in D_A, \;g\in D_{A^+}$, а $[f,g]_A(\beta)$ and
$[f,g]_A(\alpha)$ may be derived from:
$$[f,g]_A(x)=i^n\sum_{i=1}^n(-1)^{i-1}f_A^{[i-1]}(x) \overline{g_{A^+}^{[n-i]}(x)},\quad x\in I.$$

Let $\omega>0$ be a weight function from $L_{loc}(I)$,
$\omega:I\to \mathbb{R}$; the Hilbert space $L^2(I:\omega)$ is
formed as usual.

We define maximal and minimal operators as follows:
\begin{defn}\label{maxminop} Operators $T_{max}$ and $T_{min}$ are called respectively
\textit{maximal} and \textit{minimal} operators if they are
generated by $\tau(f)$ on the domains $D(T_{max})$ and
$D(T_{min})$:

\begin{equation*}\label{maximal}
D(T_{max})=\{f:I\to \mathbb{C}|\;f\in
D_A(I);\;\omega^{-1}\tau(f)\in L^2(I:\omega) \},
\end{equation*}
$$ T_{max}f=\omega^{-1}\tau(f),\;(f\in D(T_{max}));$$
\begin{equation*}\label{minimal}
D(T_{min})=\{f|\;f\in D(T_{max});\;[f,g]_{A}(b)-[f,g]_{A}(a) =
0\;(g\in D(T_{max}))\},
\end{equation*}
$$ T_{min}f=\omega^{-1}\tau(f),\;(f\in D(T_{min})),$$ where $[f,g]_{A}(b)$ and $[f,g]_{A}(a)$
are the limits (which necessarily exist) of the bilinear forms
from (\ref{lagrangeformula}), that is $\lim_{\beta\to
b}[f,g]_{A}(\beta)=[f,g]_{A}(b)$ and $ \lim_{\alpha\to
a}[f,g]_{A}(\alpha)=[f,g]_{A}(a)$.
\end{defn}

The following general theorem is known for the operators $T_{max}$
and $T_{min}$ :
\begin{thm}\label{everitttheorem}
 For the operators $T_{max}$ and $T_{min}$  and their domains the following facts are valid \label{prop}:
\\ $(a)\quad$  $D(T_{min}) \subseteq D(T_{max})$. Domains $D(T_{min})$ and $D(T_{max})$ are dense
in $L^2(I:\omega)$;\\ $(b)\quad$ The operator $T_{min}$ is closed
and
symmetric, the operator $T_{max}$ is closed in $L^2(I:\omega)$;\\
$(c)\quad$ $T_{min}^*=T_{max}$ and $T_{max}^*= T_{min}$.
\end{thm}

All self-adjoint extensions of $T_{min}$ appear to be the
contractions of $T_{max}$.

Let $\Omega$ be a finite or a countable set of indices. On
$\Omega$, we have an Everitt-Markus-Zettl multi-interval
quasi-differential system $\{I_i,\tau_i;\omega_i\}_{i\in\Omega}$.
This EMZ system generates a family of the weighted Hilbert spaces
$\{L^2(I_i:\omega_i)= L^2_i\}_{i\in\Omega}$ and families of
minimal $\{T_{min,i}\}_{i\in\Omega}$ and maximal
$\{T_{max,i}\}_{i\in\Omega}$ operators. Consider a respective
family $\{T_i\}_{i\in\Omega}$ of self-adjoint extensions.

We introduce the system Hilbert space ${\mathbf{L}^2}=\dir{
L^2_i}$ consisting of vectors $\mathbf{f}=\dir{f_i}$, such that
$f_i\in
 L^2_i$ and
\begin{equation*}
\norm{\mathbf{f}}=\sum_{i\in\Omega}\normi{f_i}=\sum_{i\in
\Omega}\int_{I_i}|f_i|^2\omega_i\,dx<\infty,
\end{equation*}
where $\normi{\cdot}$ are the norms in $ L^2_i$. In the space
${\mathbf{L}^2}$ consider the operator
$T:D(T)\subseteq{\mathbf{L}^2}\rightarrow{\mathbf{L}^2},$ defined
on the domain
\begin{equation*}\label{domain}
D(T)=\left\{\mathbf{f}\in\dir{D(T_i)}\subseteq{\mathbf{L}^2}:\sum_{i\in\Omega}\normi{T_if_i}<\infty\right\}
\end{equation*}
by $ T\mathbf{f}=\dir{T_if_i}. $
\begin{defn}
The operator $T=\oplus_{i\in\Omega}T_i$ is called a
\emph{differential vector-operator} generated by the self-adjoint
extensions $T_i$ on an \emph{EMZ} system, or simply a
vector-operator. If $\Omega$ is infinite, the vector-operator $T$
is called \emph{infinite}. The operators $T_i$ are called
\emph{coordinate} operators. For $\Omega'\subset\Omega$, the
operator $\oplus_{k\in\Omega'}T_k$ is called a
\emph{sub-vector-operator} of the vector-operator
$\oplus_{i\in\Omega}T_i$.
\end{defn}

 The following abstract preliminaries may be found, for
instance, in books \cite{reedsimon, danford}.

Fix $i\in\Omega$. For each $T_i$ there exists a unique resolution
of the identity $\iri$ and a unitary operator $U_i$, making the
isometrically isomorphic mapping of the Hilbert space $ L^2_i$
onto the space $\ltwoi$, where the operator $T_i$ is represented
as a multiplication operator. Below, we remind the structure of
the mapping $U_i$.

 We call $\phi\in L^2_i$ a \emph{cyclic vector} if for each
$z\in L^2_i$ there exists a Borel function $f$, such that
$z=f(T_i)\phi$. Generally, there is no a cyclic vector in $ L^2_i$
but there is a collection $\{\phi^k\}$ of them  in $ L^2_i$, such
that $ L^2_i=\dirk{ L^2_i(\phi^k)}$, where $ L^2_i(\phi^k)$ are
$T_i$-invariant subspaces in $ L^2_i$ generated by the cyclic
vectors $\phi^k$. That is $$
 L^2_i(\phi^k)=\overline{\{f(T_i)\phi^k\}},
$$ for a varying Borel function $f$, such that $\phi^k\in
D(f(T_i))$. There exist unitary operators $$U^k: L^2_i(\phi^k)\to
L^2(\mathbb{R},\mu^k),$$ where
$\mu^k(\Delta)={\|E^i(\Delta)\phi^k\|_i^2}$ for any Borel set
$\Delta$. In $L^2(\mathbb{R},\mu^k)$, the operator $T_i$ has the
form of multiplication by $\lambda$, i.e. $$\left(U^kT_i|_{
L^2_i(\phi^k)}{U^k}^{-1}z\right)(\lambda)=\lambda z(\lambda).$$
Then the operator $$U_i=\dirk{U^k}:\dirk{
L^2_i(\phi^k)}\to\dirk{L^2(\mathbb{R},\mu^k)}$$ makes the spectral
representation of the space $ L^2_i$ onto the space
$L^2(M_i,\mu_i)$, where $M_i$ is a union of nonintersecting copies
of the real line (\emph{a sliced union}) and $\mu_i=\sum_k\mu^k$.
That is $(U_iT_iU_i^{-1}z)(\lambda)=f(\lambda)z(\lambda),$ where
$z\in U[D(T_i)]$ and $f$ is a Borel function defined almost
everywhere according to the measure $\mu_i$.

A vector $\phi\in L^2_i$ is called \emph{maximal} relative to the
operator $T_i$, if each measure $(E^i(\cdot)x,x)_i$, $x\in L^2_i$,
is absolutely continuous relative to the measure
$(E^i(\cdot)\phi,\phi)_i$.

For each Hilbert space $ L^2_i$, there exist a unique (up to
unitary equivalence) decomposition $ L^2_i=\oplus_k
L^2_i(\varphi^k_i)$, where $\varphi^1_i$ is maximal in $ L^2_i$
relative to $T_i$, and a decreasing set of multiplicity sets
$e_k^i$, where $e_1^i$ is the whole line, such that $\oplus_k
L^2_i(\varphi^k_i)$ is equivalent with $\oplus_k
L^2(e^i_k,\mu_i)$, where the measure of the ordered representation
is defined as
$\mu_i(\cdot)=(E^i(\cdot)\varphi_i^1,\varphi_i^1)_i$. A spectral
representation of $T_i$ in $\oplus_k L^2(e^i_k,\mu_i)$ is called
the \emph{ordered representation} and it is unique, up to a
unitary equivalence. Two operators are called \emph{equivalent},
if they create the same ordered representation of their spaces.

A well-known theorem \cite[Ch. XIII, Section 5, Theorem
1]{danford}) represents the structural result for the ordered
representation of the operator $T_i$ in its abstract form. Since
the generalized eigenfunctions $W_k(x,\lambda)$ from this theorem
 are only measurable with respect to the spectral parameter $\lambda$,
the usual technique is to decompose them using an analytical basis
of solutions of the equation $(\tau_i-\lambda)\sigma=0$. At that,
frequently we do not need all the basis functions and use only a
part of them. The \emph{Defining system} $\sigma_1,\dots,\sigma_s$
is the subsystem of the solution basis such that all
$W_k(\cdot,\lambda)$ belong to its linear capsule. This treatment
leads to a very important conception of matrix Hilbert spaces.

\section{The spectral representation for the vector-operator $T$}
In this section we show, how the common spectral representation of
the vector-operator $T$ depends on the common spectral
representations of the given operators $T_i$. For this purpose, we
first prove some auxiliary results.

\begin{defn}\label{definition1}
For $i\in\Omega$, we introduce a \emph{sliced union} of sets $M_i$
(see also preliminaries) as a set $M$, containing all $M_i$ on
different copies of $\cup_{i\in\Omega} M_i$. The sets $M_i$ do not
intersect in $M$, but they can \emph{superpose}, i.e. two sets
$M_i$ and $M_j$ superpose, if their projections in the set
$\cup_{i\in\Omega} M_i$ intersect.
\end{defn}

Separate arguments show, that the following auxiliary proposition
is true.
\begin{prop}\label{proposition}
Let us have a set of measures $\mu_i$, $i\in \Omega$, defined on
nonintersecting supports. If
$$\sum_{i\in\Omega}\intinf{f(\lambda)}{\mu_i(\lambda)}<\infty,$$
 for any Borel function $f(\lambda)$, then the following equality is true:
\begin{equation*}\label{prop}
\sum_{i\in\Omega}\intinf{f(\lambda)}{\mu_i(\lambda)} =
\intinf{f(\lambda)}{\sum_{i\in\Omega}\mu_i(\lambda)}.
\end{equation*}
\end{prop}

\begin{lem}\label{irequality}
The identity resolution $\{E_\lambda\}$ of the vector-operator $T$
equals to the direct sum of the coordinate identity resolutions
$\{\iri\}$, that is:
\begin{equation*}
\{E_\lambda\} = \dir{\{\iri\}}
\end{equation*}
\end{lem}

\begin{proof}
A vector $\mathbf{x}$ belongs to $D(T)$ if and only if
$$\norm{T\mathbf{x}}=\sum_{i\in\Omega}\normi{T_ix_i}=\sum_{i\in\Omega}\intinf{\lambda^2}{\normi{\iri
x_i}}<\infty.$$  Then, using Proposition $\ref{proposition}$ we
find out that:
\begin{equation*}\label{walkingsum}
\sum_{i\in\Omega}\intinf{\lambda^2}{\normi{\iri
x_i}}=\intinf{\lambda^2}{\sum_{i\in\Omega}\normi{\iri x_i}}.
\end{equation*}
This means, that $\mathbf{x}\in D(T)$, if and only if
$$\intinf{\lambda^2}{\sum_{i\in\Omega}\normi{\iri x_i}}<\infty$$ and
$$\norm{T\mathbf{x}}=\intinf{\lambda^2}{\sum_{i\in\Omega}\normi{\iri
x_i}}.$$ Using the  uniqueness property of an identity resolution,
the last two equations show that $\dir{\{\iri\}}$ is the identity
resolution of the vector-operator $T$. That is, according to our
notations $\{E_\lambda\}=\dir{\{\iri\}}$. The lemma is proved.
\end{proof}
\begin{lem}\label{Borel function}
For any Borel function $f$ and any vector $\mathbf{x}\in D(f(T))$,
the following equality holds:
$f(T)\mathbf{x}=[\dir{f(T_i)}]\mathbf{x}$.
\end{lem}
\begin{proof}
Let $\mathbf{x}\in D(f(T))$. Then, paying attention to Proposition
$\ref{proposition}$ and Lemma $\ref{irequality}$, for any
$\mathbf{y}\in \mathbf{L}^2$, we obtain:
\begin{multline*}
(f(T)\mathbf{x},\mathbf{y})=\intinf{f(\lambda)}{(E_\lambda
\mathbf{x},\mathbf{y})}=
\intinf{f(\lambda)}{\sum_{i\in\Omega}(E^i_\lambda x_i,y_i)_i}=\\
=\sum_{i\in\Omega}\intinf{f(\lambda)}{(E^i_\lambda
x_i,y_i)_i}=\sum_{i\in\Omega}(f(T_i)x_i,y_i)_i=([\dir{f(T)}]\mathbf{x},\mathbf{y}).
\end{multline*}
Since $\mathbf{y}$ is arbitrary, we have
$f(T)\mathbf{x}=[\dir{f(T_i)}]\mathbf{x}$. This completes the
proof of the lemma.
\end{proof}
For $z_i\in  L^2_i$, $i\in\Omega$, define
$\overline{\mathbf{z_i}}=\{0,...,0,z_i,0,...,0\}\in{\mathbf{L}^2},$
where $z_i$ is on the $i$-th place.

For each $i\in\Omega$, let $\epsilon(T_i)$ denote the
\emph{subspectrum} of the operator $T_i$, i.e. the set where the
spectral measures of $T_i$ are concentrated. Note, that
$\overline{\epsilon(T_i)}=\sigma(T_i)$. For instance, the
subspectrum of an operator having the complete system of
eigenfunctions with eigenvalues being the rational numbers of
$[0,1]$ equals to $\mathbb{Q}\cap [0,1]$; the subspectrum of an
operator having the continuous spectrum [0,1] is assumed to equal
to (0,1) without loss of generality.

 Consider a projecting mapping $P:M\to \cup_{i\in\Omega} M_i$ (see
Definition \ref{definition1}), such that
$P(\epsilon(T_i))=\epsilon(T_i)$.

\begin{defn}\label{ak}
Let $\Omega=\cup_{k=1}^K A_k$, $A_k\cap A_s=\emptyset$ for $k\neq
s$ and
\begin{multline*}
A_k=\{s\in\Omega: \forall s,l\in A_k, s\neq l,
P(\epsilon(T_s))\cap P(\epsilon(T_l))=B_{sl},\\ \mbox{where}\,
\norm{E^t(B_{sl})\varphi_t}_t=0 \,\,\mbox{for any cyclic}\,\,
\varphi_t \in L^2_t,\,t=s,l\}.
\end{multline*}
 From all such divisions of $\Omega$ we choose and fix the
one, which contains the minimal number of $A_k$. In case all the
coordinate spectra $\sigma(T_i)$ are simple, we define the number
$\Lambda=\min\{K\}$ as the \emph{spectral index} of the
vector-operator $T$.
\end{defn}

\begin{thm}\label{spectral index}
Let each $T_i$ have a cyclic vector $a_i$ in $ L^2_i$. Then the
vector-operator $T$ has $\Lambda$ cyclic vectors
$\{\mathbf{a}_k\}_{k=1}^\Lambda$, having the form $
\mathbf{a}_k=\sum_{i\in A_k}\overline{\mathbf{a_i}}. $
\end{thm}
\begin{proof}
First we consider the case of two coordinate operators. Let
$s,l\in \Omega$. Then, in order to obtain one cyclic vector in
$L^2_s\oplus L^2_l$ having the form $a_s\oplus a_l$, for any
$\mathbf{x}=x_s\oplus x_l\in L^2_s\oplus L^2_l$ we have to find a
Borel function $f$, such that
$$\mathbf{x}=f(T_s\oplus T_l)[a_s\oplus a_l].$$ From Lemma
$\ref{Borel function}$ it follows that $$\mathbf{x}=[f(T_s)\oplus
f(T_l)][a_s\oplus a_l].$$ On the other hand we must obtain each
space $L^2_p\, (p=s,l)$ by closing the set $\{f_p(T_p)a_p\}$,
letting $f_p$ vary over all the Borel functions such that $a_p\in
D(f_p(T_p))$. If $s,l\in A_k$, then supposing that $f=f_p$ on
$P(\epsilon(T_p))$, we obtain the required function $f$, since any
functions in the isomorphic space $L^2$ are considered equal on
the set of measure zero. Hence, it is clear that for all $i\in
A_k$, we may build a single cyclic vector of the form
$$\mathbf{a}_k=\oplus_{i\in A_k}a_i = \sum_{i\in A_k}\overline{\mathbf{a_i}},$$
using the process described above, each time operating with a pair
of operators.

We recall, that we have the minimal number of $A_k$. Consider the
Hilbert space
\begin{equation}\label{2hs}
[\oplus_{i\in A_k} L^2_i]\oplus[\oplus_{j\in A_q}L^2_j], k\neq q.
\end{equation}
We know, that then $$[\cup_{i\in A_k}P(\epsilon(T_i))]\cap
[\cup_{j\in A_q}P(\epsilon(T_j))] = B_{kq}$$ has a non-zero
spectral measure. From the reasonings described in the beginning
of this proof we see, that for joining the cyclic vectors
$\mathbf{a}_k=\oplus_{i\in A_k}a_i$ and $\mathbf{a}_q=\oplus_{j\in
A_q}a_j$ into the one
$$\mathbf{a}_k+\mathbf{a}_q=\sum_{i\in A_k}\overline{\mathbf{a_i}} +
\sum_{j\in A_q}\overline{\mathbf{a_j}},$$ we would have to derive
the Hilbert space ($\ref{2hs}$) by closing the set
$$\{f_k(\oplus_{i\in
A_k}T_i)\mathbf{a}_k\}\oplus\{f_q(\oplus_{j\in
A_q}T_j)\mathbf{a}_q\},$$ with varying the Borel functions $f_k$
and $f_q$, which coincide on $B_{kq}$. This is not possible, since
the set of such functions is not dense in the isomorphic space
$L^2$ (the isomorphism is understood as in the spectral
representation of the space ($\ref{2hs}$)). Hence, we have
obtained $\Lambda$ cyclic vectors
$$\mathbf{a}_k=\sum_{i\in A_k}\overline{\mathbf{a_i}}\in
{\mathbf{L}^2},\,k=\overline{1,\Lambda}$$ and have proved the
theorem.
\end{proof}

\begin{cor}
Let each $T_i$ have a single cyclic vector.  Then \\
1. $\Lambda=1$ if and only if the coordinate operators $T_i, i\in
\Omega$, have almost everywhere (relative to the spectral measure)
pairwise non-superposing subspectra. \\2. a)
$\mathrm{card}(\Omega)<\aleph_0$. $\Lambda=\mathrm{card}(\Omega)$,
if and only if all the coordinate operators $T_i$ have pairwise
superposing subspectra; b) $\mathrm{card}(\Omega)=\aleph_0$.
$\Lambda=\infty$, if and only if $T$ has an infinite
sub-vector-operator, the coordinate operators of which have
pairwise superposing subspectra.
\end{cor}
\begin{proof}
The proof directly follows from the reasonings of the proof of
Theorem $\ref{spectral index}$.
\end{proof}

In the next section we will rigorously show what the spectral
multiplicity of a vector-operator is. Nevertheless, this notation
is intuitively clear. Running ahead, let us present here an
example, which will show the difference between the spectral index
and the spectral multiplicity of the vector-operator $T$.

\textbf{Example 1.} We have a three-interval EMZ system
$\{I_i,\tau_i\}_{i=1}^3$ (a kinetic energy, a mirror kinetic
energy, an impulse):
\begin{multline*}I_1=[0,+\infty),\,\,\,\tau_1=-\left(\frac
d{dt}\right)^2,\;\;\\D(T_1)=\{f\in D(T_{max,1}): f(0)+kf'(0)=0,
-\infty<k\leqslant\infty\};\end{multline*}
\begin{multline*}I_2=[0,+\infty),\,\,\,\tau_2=\left(\frac
d{dt}\right)^2,\;\;\\D(T_2)=\{f\in D(T_{max,2}): f(0)+sf'(0)=0,
-\infty<s\leqslant\infty\};\end{multline*}
\begin{multline*}I_3=[0,1],\,\,\,\tau_3=\frac 1i\frac
d{dt},\;\;D(T_3)=\{f\in D(T_{max,3}):
f(0)=e^{i\alpha}f(1),\,\alpha\in [0,2\pi]\}.\end{multline*} If
$k,s\in (-\infty,0]\cup\{+\infty\}$ then
$$\epsilon(T_1)=(0,+\infty),\;\;\epsilon(T_2)=(-\infty,0),\;\;\epsilon(T_3)=\bigcup_{n=-\infty}^\infty (2\pi n -\alpha).$$
For this system we have: $\{1,2,3\}=\cup_{k=1}^2 A_k$ and
$A_1=\{1,2\},\,A_2=\{3\}$. Thus, here the spectral index does not
coincide with the spectral multiplicity (which equals to 1) and
equals to 2.

The case $0<k,s<+\infty$ leads to the following
$$
\epsilon(T_1)=\left\{-\frac 1 {k^2}\right\}\cup(0,+\infty),
\;\;\epsilon(T_2)= (-\infty,0)\cup \left\{\frac 1
{s^2}\right\},\;\;\epsilon(T_3)=\bigcup_{n=-\infty}^\infty (2\pi n
-\alpha). $$
 If $$\alpha\not\in \left[\bigcup_{n=-\infty}^\infty
\left(2\pi n +\frac1{k^2}\right)\right]\bigcup\,
\left[\bigcup_{n=-\infty}^\infty \left(2\pi n
-\frac1{s^2}\right)\right],$$ we have
$A_1=\{1\},\,A_2=\{2\},\,A_3=\{3\}.$ That is $\Lambda=3$ but
$\oplus_{i=1}^3T_i$ has a simple spectrum.

\textbf{Example 2.} Let us have a vector-operator
$\oplus_{i=1}^3T_i$ with $$\epsilon(T_1)=\bigcup_{n\in\mathbb{Z},
n\geqslant 0} n,\,\,\,\epsilon(T_2)=\bigcup_{n\in\mathbb{Z},
n\leqslant 0} n,\,\,\,\epsilon(T_3)=\bigcup_{n\in\mathbb{Z}, n\neq
0} n.$$ Spectral index equals to 3 but spectral multiplicity
equals to 2.

\begin{defn}\label{distorted def}
A vector-operator $T=\oplus_{i\in\Omega}T_i$ with simple
coordinate spectra $\sigma(T_i)$ is called \emph{distorted} if its
spectral index does not equal its spectral multiplicity.
\end{defn}

Note that in the above Example 1, it is possible to unite the
cyclic vectors into one just taking their direct sum (as it is
shown in the proof of Theorem \ref{spectral index}). But
nevertheless, it is convenient to consider such operators as
distorted satisfying Definition \ref{distorted def}. The distorted
vector-operator from Example 2 is 'completely' distorted and it is
not possible to unite the coordinate cyclic vectors into a cyclic
direct sum.

With some loss of technical value but more clearly for
applications, Theorem \ref{spectral index} may be reformulated as
\begin{cor}
Let each $T_i$ have a simple spectrum. Then undistorted
vector-operator $T$ has $\Lambda$-multiple spectrum.
\end{cor}
Let us pass to the general case when each operator $T_i$ has $m_i$
cyclic vectors. There exists a decomposition
\begin{equation*}\label{decomposition}
T=\dir{T_i}=\dir{\oplus_{k=i}^{m_i}T_i^k}=\oplus_{s}T_s,
\end{equation*}
where each $T_s$ has a single cyclic vector. For the
vector-operator $T$ decomposed as above, we apply Theorem
$\ref{spectral index}$ and find the spectral index $\Lambda$. It
is clear, that in this case for the spectral index there exists
the estimate
\begin{equation}
\Lambda\geqslant \max\{m_i\}.
\end{equation}
As it has been stated in the preliminaries, for each operator
$T_i$ there exists the unitary operator $U_i$, such that
$U_i:L^2_i\to L^2(M_i,\mu_i)$. Hence
$$\dir{U_i}:\dir{L^2_i}\to\dir{L^2(M_i,\mu_i)}.$$
Or, in the general case (i.e. when there are $T_i$ with more then
one cyclic vector),
$$\dir{U_i}:\dir{\oplus_{k=1}^{m_i}L^2_{i,k}}\to
\dir{\oplus_{k=1}^{m_i} L^2(\mathbb{R},\mu^k_i)}.$$ From Theorem
$\ref{spectral index}$ it follows that there exists a unitary
operator
\begin{equation}\label{uov}
V:\dir{\oplus_{k=1}^{m_i} L^2(\mathbb{R},\mu^k_i)}= \oplus_s
L^2(\mathbb{R},\mu_s)\to\oplus_{q=1}^\Lambda
L^2\left(\mathbb{R},\sum_{j\in A_q}\mu_j\right).
\end{equation}
This means that for any vector-operator $T$ there exists the
unitary operator $V\dir{U_i}$, which represents the space
${\mathbf{L}^2}$ on the space $L_2(N,\mu)$:
\begin{equation*}
V\dir{U_i}:{\mathbf{L}^2}\to L^2(N,\mu),
\end{equation*}
where $N$ is the sliced union of $\Lambda$ copies of $\mathbb{R}$
and $$\mu=\sum_{q=1}^\Lambda\sum_{j\in A_q}\mu_j,$$ according to
the symbols in $(\ref{uov})$. We finally obtain
\begin{thm}\label{spectral representation}
Let the vector-operator $T=\dir T_i$ be undistorted and let the
unitary operator $V$ be defined as in (\ref{uov}). If unitary
operators $U_i$ give spectral representations of the Hilbert
spaces $L^2_i$ on the spaces $\ltwoi$, then the unitary operator
$$W=V\dir{U_i}$$ gives a spectral representation of the space ${\mathbf{L}^2}$
on the space $L^2(N,\mu)$.
\end{thm}

Directly from the definition of a distorted vector-operator, it
follows that only for undistorted vector-operators, the transform
$V$ does reduce the quantity of cyclic vectors to the minimal
possible. Note, that distorted differential vector-operators
appear to be frequent objects if vector-operators are considered
on a set of closed bounded intervals, and on the other hand quite
rare, if coordinate operators have continuous spectra. For them
Theorem \ref{spectral representation} is not efficient and needs
to be strengthened. Such strengthening is the construction of an
ordered representation for arbitrary (distorted or not)
differential vector-operators, the process which seems to be
essential for further development of spectral theory of
vector-operators.

\section{The ordered spectral representation for the vector-operator $T$}

\begin{thm}\label{ordered representation}
If $\theta_i$ and $\{e_n^i\}_{n=1}^{m_i}$ are measures  and
multiplicity sets of ordered representations for coordinate
operators $T_i$, $i\in \Omega$, then there exist processes $Pr_1$
and $Pr_2$, such that the measure $$\theta =
Pr_1(\{\theta_i\}_{i\in\Omega})$$ is the measure of an ordered
representation and the sets $$s_n =
Pr_2\left(\{e_k^i\}_{i\in\Omega;\, k=\overline{1,m_i}}\right)$$
are the canonical multiplicity sets of the ordered representation
of the operator $T$. Thus, the unitary representation of the space
${\mathbf{L}^2}$ on the space $\oplus_n L^2(s_n,\theta)$ is the
ordered representation and it is unique up to unitary equivalence.
\end{thm}
\begin{proof}
We divide the proof into units for convenience. Units \textbf{(A)}
and \textbf{(B)} represent the process, which we call 'the process
of division on subspectra'.

\textbf{(A)} Let $a_i$ be maximal vectors relative to the
operators $T_i$ in $L^2_i$. We want to find a maximal vector
relative to the vector-operator $T$. We know, that the vector
$\dir a_i$ does not give a single measure, if a set
$P(\epsilon(T_i))\cap P(\epsilon(T_j))$ has a non-zero spectral
measure for $i\neq j$. Consider restrictions $T_i|_{
L^2_i(a_i)}=T_i'$. Since all the operators $T_i'$ have single
cyclic vectors $a_i$, we can divide $\Omega$ into $A_k$,
$k=\overline{1,\Lambda}$ (see definition \ref{ak}) and apply
Theorem \ref{spectral index} for the operator $\dir T_i'$. Thus,
we have derived $\Lambda$ vectors $\mathbf{a}^k = \oplus_{j\in
A_k}a_j$, which are maximal in the respective spaces
${\mathbf{L}^2}(\mathbf{a}^k)=\oplus_{j\in A_k}L^2_j(a_j)$.
Indeed, this is obvious for the case
$\mathrm{card}(A_k)<\aleph_0$. For the infinite case, if arbitrary
$\mathbf{y}=\oplus_{j\in A_k}y_j\in {\mathbf{L}^2}(\mathbf{a}^k)$
and if

\begin{equation}\label{conv}\left([\oplus_{j\in
A_k}E^j](\cdot)\mathbf{a}^k,\mathbf{a}^k\right)=\sum_{j\in
A_k}\left(E^j(\cdot)a_j,a_j\right)_j=0, \end{equation}
 then from the
maximality of the vectors $a_j$ for all $j\in A_k$, and since
$P(\epsilon(T_j'))\cap P(\epsilon(T_k'))$ has zero spectral
measures for $j\neq k$, we obtain
$$\sum_{j\in A_k}\left(E^j(\cdot)y_j,y_j\right)_j=\left([\oplus_{j\in
A_k}E^j](\cdot)\mathbf{y},\mathbf{y}\right)=0,$$ which follows
from the convergence to zero of the series with the positive
maximal elements (\ref{conv}). Thus, in particular, we have
constructed a maximal vector in ${\mathbf{L}^2}$ for the case
$\Lambda=1$.

\textbf{(B)} Let now $1<\Lambda<\infty$. Define $T^k=\oplus_{j\in
A_k}T_j'$. For any two operators $T^k$ and $T^s$, $k\neq s$, let
us introduce the sets $\epsilon_{k,s}=P(\epsilon(T^k))\cap
P(\epsilon(T^s))\,\,\,\mbox{and}\,\,\,\epsilon_{k}=P(\epsilon(T^k))\setminus
\epsilon_{k,s}.$ There exist unitary representations
$U^k:{\mathbf{L}^2}(\mathbf{a}^k)\rightarrow
L^2(\mathbb{R},\mu_{\mathbf{a}^k})$ (see formula (\ref{uov})
supposing there  $\Lambda=1$). Consider measures $\mu_k$ and
$\mu_{k,s}$, defined as
$$\mu_{k,s}(e)=\mu_{\mathbf{a}^k}(e\cap\epsilon_{k,s})\,\,\,\mbox{and}\,\,\,\mu_k(e)=\mu_{\mathbf{a}^k}(e\cap\epsilon_k),$$
for any measurable set $e$. For any operator $T^k$ (with respect
to $T^s$), measures $\mu_k$ and $\mu_{k,s}$ are mutually singular
and $\mu_k+\mu_{k,s}=\mu_{\mathbf{a}^k}$; therefore
$$L^2(\mathbb{R},\mu_{\mathbf{a}^k})=L^2(\mathbb{R},\mu_k)\oplus
L^2(\mathbb{R},\mu_{k,s}).$$ This means that (according to our
designations):
$${U^k}^{-1}:L^2(\mathbb{R},\mu_{\mathbf{a}^k})\longrightarrow
\mathbf{L}^2(\mathbf{a}^k_k)\oplus
\mathbf{L}^2(\mathbf{a}^k_{k,s})$$ and
$\mathbf{a}^k=\mathbf{a}^k_k\oplus \mathbf{a}^k_{k,s}$, where
$\mathbf{a}^k_k$ and $\mathbf{a}^k_{k,s}$ form the measures
$\mu_k$ and $\mu_{k,s}$ respectively. Define also as $\max\{w,
\psi\}$ the vector, which is maximal of the two vectors in the
brackets (Note, that this designation is valid only for vectors,
considered on the same set. In order not to complicate the
investigation we assume here that any two vectors are comparable
in this sense. In order to achieve this, it is enough to decompose
each coordinate operator $T_i$ into the direct sum $T_i^{pp}\oplus
T_i^{cont}$, where the operators have respectively pure point and
continuous spectra. Then after redesignation we obtain the
equivalent vector-operator to the initial vector-operator $\oplus
T_i$).

Consider first two operators $T^1$ and $T^2$. It is clear, that
the vector
$$\mathbf{a}^{1\oplus 2}=\mathbf{a}^1_1\oplus \mathbf{a}_2^2\oplus\max\left\{\mathbf{a}^1_{1,2},\mathbf{a}^2_{2,1}\right\}$$ is
maximal in
$\mathbf{L}^2(\mathbf{a}^1)\oplus\mathbf{L}^2(\mathbf{a}^2)$.
Note, that $\mathbf{a}^1_1$ and $\mathbf{a}_2^2$ and they both may
equal zero. The maximal vector in
$\mathbf{L}^2(\mathbf{a}^1)\oplus\mathbf{L}^2(\mathbf{a}^2)\oplus\mathbf{L}^2(\mathbf{a}^3)\
$ will have the form:
$$\mathbf{a}^{1\oplus 2\oplus 3}=\mathbf{a}^{1\oplus 2}_{1\oplus
2}\oplus \mathbf{a}_3^3\oplus\max\left\{\mathbf{a}^{1\oplus
2}_{{1\oplus 2},3},\mathbf{a}^3_{3,{1\oplus 2}}\right\}.$$
Continuing this process, we obtain a maximal vector in the main
space ${\mathbf{L}^2}$:
\begin{equation}\label{max-vec}
\mathbf{a}^{1\oplus \cdots \oplus
\Lambda}=\mathbf{a}^{1\oplus\cdots \oplus\Lambda-1}_{1\oplus\cdots
\oplus\Lambda-1}\oplus
\mathbf{a}_\Lambda^\Lambda\oplus\max\left\{\mathbf{a}^{1\oplus\cdots
\oplus\Lambda-1}_{{1\oplus\cdots
\oplus\Lambda-1},\Lambda},\mathbf{a}^\Lambda_{\Lambda,{1\oplus\cdots
\oplus\Lambda-1}}\right\}.
\end{equation}

Formula \ref{max-vec} may be simplified, if we divide the measures
$\mu_{\mathbf{a}^k}$ into continuous and pure point components,
that is
$\mu_{\mathbf{a}^k}=\mu^{cont}_{\mathbf{a}^k}+\mu_{\mathbf{a}^k}^{pp}$.
Then $\mathbf{a}^k=\mathbf{a}^{k,cont}\oplus \mathbf{a}^{k,pp}$.
Relatively to any operator $T^s$, $k\neq s$, we have
$$\mathbf{a}^{k,cont}=\mathbf{a}_k^{k,cont}\oplus \mathbf{a}_{k,s}^{k,cont}\,\,\mbox{and}\,\,\mathbf{a}^{k,pp}=\mathbf{a}_{k}^{k,pp}\oplus \mathbf{a}_{k,s}^{k,pp}.$$
Now we can repeat the process described above in \textbf{(B)},
separately for the continuous and the pure point parts. Since
measures with the same null set may be considered equivalent, we
have
$$\max\{w^{cont},
\psi^{cont}\}=\,\mbox{either}\,w^{cont}\,\mbox{or}\,\psi^{cont},$$
$$\max\{w^{pp},
\psi^{pp}\}=\,\mbox{either}\,w^{pp}\,\mbox{or}\,\psi^{pp},$$ for
any two vectors $w$ and $\psi$. Thus we obtain $$
\mathbf{a}^{1\oplus \cdots \oplus
\Lambda,cont}=\mathbf{a}^{1,cont}\oplus\left[\oplus_{j=2}^\Lambda
\mathbf{a}_j^{j,cont}\right]. $$ Similarly, $$ \mathbf{a}^{1\oplus
\cdots \oplus
\Lambda,pp}=\mathbf{a}^{1,pp}\oplus\left[\oplus_{j=2}^\Lambda
\mathbf{a}_j^{j,pp}\right]. $$ Since $\max\{w^{cont},
\psi^{pp}\}=\psi^{pp},$ we finally derive
\begin{equation*}\label{cont-max}
\mathbf{a}^{1\oplus \cdots \oplus \Lambda}=\mathbf{a}^{1\oplus
\cdots \oplus \Lambda,pp}\oplus \mathbf{a}^{1\oplus \cdots \oplus
\Lambda,cont}_{1\oplus \cdots \oplus \Lambda}.
\end{equation*}

Let $\Lambda=\infty$. We obtain $\mathbf{a}^{1\oplus \cdots \oplus
\Lambda}$ as a vector which satisfies the following equality:
\begin{equation}\label{convmaxv}\left\|[\dir{E^i(\cdot)}]\,\mathbf{a}^{1\oplus \cdots \oplus
\Lambda}\right\|^2=\lim_{L\to\infty}
\left\|[\oplus_{j=1}^L{E^j(\cdot)}]\,\mathbf{a}^{1\oplus \cdots
\oplus L}\right\|^2,\end{equation} since the limit on the right
side exists. Indeed $\lim_{L\to\infty}
\left\|[\oplus_{j=1}^L{E^j(\cdot)}]\,\mathbf{a}^{1\oplus \cdots
\oplus L}\right\|^2$ can be rewritten as
$\sum_{j=1}^\infty\left\|E^j(\cdot)\,\widehat{a_j}\right\|_j^2$,
where $\widehat{a_j}$ are the restricted $a_j$. Noticing that
$$\sum_{j=1}^\infty\left\|{E^j(\cdot)}\,\widehat{a_j}\right\|_j^2\leqslant
\sum_{j=1}^\infty \left\|{E^i(\cdot)}\,a_j\right\|_j^2<\infty,$$
we prove the convergence (without loss of generality, the vectors
$a_i$ can be always chosen such, that
$\sum_{i=1}^\infty\|a_i\|_i^2<\infty$).

\textbf{(C)} The next step is to build the measure of the ordered
representation for the vector-operator. From Lemma
\ref{irequality} and the reasonings above, it follows that such a
measure will be
\begin{equation*}\label{measure}
\theta(\cdot)=\left([\dir{E^i(\cdot)}]\,\mathbf{a}^{1\oplus \cdots
\oplus \Lambda},\mathbf{a}^{1\oplus \cdots \oplus \Lambda}\right).
\end{equation*}

Thus we have constructed the process $Pr_1$.

\textbf{(D)} The final step is to construct the canonical
multiplicity sets $s_n$ of the vector-operator; $s_1$ is the whole
line; $s_2$ must contain all the spectrum the multiplicity of
which exceeds or equals to 2. For this purpose, we are primarily
to unite all $e_2^i$. But, nevertheless, $\cup_i e_2^i$ will not
include all the sets of multiplicity $\geqslant 2$, since we know
that if $P(e_1^i\setminus e_2^i)\cap P(e_1^j\setminus e_2^j)$ has
a non-zero spectral measure, all the intersections of this sort
will represent the multiplicity 2 and should be included into
$s_2$ (since then it is not possible to construct a single cyclic
vector). That is $s_2=\left(\cup_i P(e_2^i)\right)\cup \left(\cup
\cap(P(e^i_1\backslash e^i_2)\right)$. Using this idea and the
fact that an infinite intersection of measurable sets is a
measurable set, by induction we may finally build $s_n$:
\begin{equation}\label{mult-set}
    s_n= \left [
\bigcup_iP(e_n^i)\right ] \bigcup \left [\bigcup_{\sum
m_i\geqslant n}\bigcap P\left(e^i_{m_i}\backslash
e^i_{m_i+1}\right)\right ].
\end{equation}

We have constructed the process $Pr_2$.
\end{proof}

The constructed measure and the multiplicity sets induce the
ordered representation. It is known that such a representation is
unique up to unitary equivalence.

Let us return to the Example 2. For the distorted vector-operator
$T_1\oplus T_2\oplus T_3$, two spectral measures will be
constructed on vectors $\mathbf{a}^{1\oplus 2\oplus 3}$ and
$$\min\{a_{1,2}^1,a_{2,1}^2\}\oplus
\min\{a_{2,3}^2,a_{3,2}^3\}\oplus \min\{a_{3,1}^3,a_{1,3}^1\},$$
where the sense of the minimums is clear.

Here the term 'distorted vector-operator' is clearly explained by
the form of its cyclic vectors. The multiplicity set $e_2$ will be
$$[P(\epsilon(T_1))\cap P(\epsilon(T_2))]\cup
[P(\epsilon(T_1))\cap P(\epsilon(T_3))]\cup [P(\epsilon(T_2))\cap
P(\epsilon(T_3))].$$

Using the obtained spectral representation we can construct
equivalence classes in families of self-adjoint operators:
\begin{defn}Two families of self-adjoint extensions $\{T_i\}_{i=1}^N$ and
$\{S_j\}_{j=1}^L$ are called equivalent, if respective
vector-operators $\oplus_{i=1}^N T_i$ and $\oplus_{j=1}^L S_j$ are
equivalent.
\end{defn}
Note, that if two families $\{T_i\}_{i=1}^N$ and $\{S_j\}_{j=1}^L$
are equivalent, it is not necessarily the case that $N=L$ and
$T_i$ is equivalent with $S_i$.

\section{The case of  coordinate differential operators}

Up to now, we have not used the structure of the coordinate
operators as differential operators. In this section we make
precise the ordered representation obtained in the previous
section.

Let $I=\bigvee_{i\in\Omega}I_i$ denote the sliced union of
intervals $I_i$.  Similarly, $I^k =\bigvee_{j\in A_k}I_j$. If
$x_i$ are variables on $I_i$, then $\vee x_i$ will designate a
variable either on $I$ or $I^k$ depending on the context. This
notation shows, that a vector-function $$z=\{z_1(x_1),\dots,
z_n(x_n),\dots\}$$ on $I$ or $I^k$ may be written as $z(\vee
x_i)$. In particular, we may also write $\mathbf{z}(\vee x_i)$
instead of $\mathbf{z}=\dir z_i$.

Let us introduce the space $\dir L^{\infty}(I_i^{n})$. Here,
$\mathbf{z}(\vee x_i) \in \dir L^{\infty}(I_i^{n})$ means that
$$\sup_{i\in\Omega} \left\{\mathrm{ess}
\sup_{x_i\in I_i^{n}}
|z_i(x_i)\chi_{I_i^{n}}(x_i)|\right\}<\infty,
$$
where for each $i$, families $\{ I_i^{n}\}_{n=1}^\infty$ represent
compact subintervals of $I_i$, such that $\cup_{n=1}^\infty
I_i^{n} = I_i$ and $\chi$ is the characteristic function. In
\cite[Lemma 2.1]{ashuroveveritt}, it was shown that $\dir
L^{\infty}(I_i^{n})=(\dir{L^1(I_i^n)})^*$, where
 the space of Lebesgue-integrable
vector-functions $\dir{L^1(I_i^n)}$ is defined analogously to
$\mathbf{L}^2$.

We also need to introduce a symbolic integral $\int_{\bigvee J_i}
f(\vee x_i)\, d(\vee x_i)$ defined by:
$$\int_{\bigvee J_i} f(\vee x_i)\, d(\vee x_i) = \oplus_i \int_{J_i} f_i(x_i)\,dx_i,$$
where $f(\vee x_i)$ is understood to be measurable relative to $
d(\vee x_i)$, if and only if $f_i(x_i)$ are measurable relative to
Lebesgue measures $dx_i$. Then $$\int_{\bigvee J_i} f(\vee x_i)\,
d(\vee x_i)<\infty$$ if and only if $\sup_i \int_{J_i}
f_i(x_i)\,dx_i <\infty$.

\begin{thm}\label{eigexpthm}
Let $T$ be a self-adjoint vector-operator, generated by an EMZ
system $\{I_i,\tau_i\}_{i\in\Omega}$. Let $U$ be an ordered
representation of the space
$\mathbf{L}^2=\oplus_{i\in\Omega}L^2(I_i)$ relative to $T$ with
the measure $\theta$ and the multiplicity sets $s_k,\,
k=\overline{1,m}$. Then there exist kernels $\Theta_k(\vee x_i,
\lambda)$, measurable relative to $d(\vee x_i)\times \theta$, such
that $\Theta_k(\vee x_i, \lambda)=0$ for $\lambda\in
\mathbb{R}\setminus s_k$ and
$(\oplus_{i\in\Omega}\tau_i-\lambda)\Theta_k(\vee x_i, \lambda)=0$
for each fixed $\lambda$. Moreover,
\begin{equation}\label{essbound1}
\int_\Delta |\Theta_k(\vee x_i, \lambda)|^2\,d\theta(\lambda)\in
\oplus_{i\in \Omega} L^{\infty}(I_i^{n})\,\,\, \forall n\in
\mathbb{N}.
\end{equation}
\begin{equation}\label{untrans1}
    (U \mathbf{w})^k(\lambda)=\lim_{n\to\infty}\int_{I^{n}}\,\mathbf{w}(\vee x_i)\,\overline{\Theta_k(\vee x_i, \lambda)}\,d(\vee
    x_i),
    \,\,\,\,\, \mathbf{w}\in \mathbf{L}^2,
\end{equation}
where the limit exists in $L^2(s_k,\theta)$. The kernels
$\{\Theta_k(\vee x_i, \lambda)\}_{k=1}^n$, $n\leqslant m$, are
 linearly independent as vector-functions of
the first variable almost everywhere relative to the measure
$\theta$ on $s_n$.
\end{thm}
\begin{proof} Fix i. If $\theta_i$ and $\{e^i_p\}_{p=1}^{m_i}$ are respectively the measure and
the multiplicity sets of an ordered representation for $T_i$, then
there exists the decomposition
$L^2_i=\oplus_{p=1}^{m_i}L^2(e_p^i,\theta_i)$, which implies
$T_i=\oplus_{p=1}^{m_i}T_i^p$ and $L^2(e_p^i,\theta_i)$ are
$T_i^p$-invariant. For vector-operator
$(\oplus_{i\in\Omega}\oplus_{p=1}^{m_i}T_i^p)\to\mbox{redesignate}\to
\oplus_sT_s,\, s=\{i,p\}\in\Omega_1$, we may write
$\Omega_1=\cup_{k=1}^\Lambda A_k$.

Let us separate the proof into units for convenience.

\textbf{(A)} For each $T_j$, $j\in A_k\subset\Omega_1$ and
$k=\overline{1,\Lambda}$, there exists a single cyclic vector
$a_j\in L^2_j$ and \cite[XII.3, Lemma 9 and XIII.5, Theorem
1(I)]{danford} a function $W_j(x_j,\lambda)$ defined on $I_j\times
e_j$ (note, that for a fixed $i\in\Omega$, $I_j=I_i$ for all
$p=\overline{1,m_i}$) and measurable relative to $dx_j\times
\mu_{a_j}$, such that $W_j(x_j,\lambda)=0,\,\lambda\in
\mathbb{R}\setminus e_j$ and for any bounded $\Delta\subset e_j$:
$$\int_\Delta |W_j(x_j,\lambda)|^2\,d\mu_{a_j}(\lambda)\in L^\infty(I_j^n),\,\,\, n\in\mathbb{N}.$$
Also
\begin{equation}\label{ran1}
    \left(E^j(\Delta)F_j(T_j)a_j\right)(x_j)=\int_\Delta
    W_j(x_j,\lambda)F_j(\lambda)\,d\mu_{a_j}(\lambda),
\end{equation}
for any $F_j\in L^2(e_j,\mu_{a_j})$. On $I^k=\bigvee_{j\in
A_k}I_j$, we construct the vector-function $$W^k(\vee x_j,
\lambda)=\{W_1(x_1,\lambda),\dots, W_n(x_n,\lambda),\dots\},$$
which is obviously measurable relative to $d(\vee x_j)\times
\sum\mu_{a_j}$. Since $W_j(\cdot,\lambda)\in
L^2(\Delta,\mu_{a_j})$, then substituting
$\overline{W_j(\cdot,\lambda)}=\overline{W_j(\lambda)}$ in
(\ref{ran1}) in place of $F_j$, we obtain
$$ \left(E^j(\Delta)\overline{W_j(T_j)}a_j\right)(x_j)=\int_\Delta |W_j(x_j,\lambda)|^2\,d\mu_{a_j}(\lambda).$$
Remembering, that $P(\epsilon(T_s))\cap P(\epsilon(T_j))$ has zero
measure, for $s\neq j$ and $s,j\in A_k$, we obtain
\begin{equation*}\label{ran}
\left(\left[\oplus_{j\in
A_k}E^j\right](\Delta)\overline{W^k(\oplus_{j\in
A_k}T_j)}\mathbf{a}^k\right)(\vee x_j)=\int_\Delta |W^k(\vee
x_j,\lambda)|^2\,d\mu_{\mathbf{a}^k}(\lambda),
\end{equation*}
where $\mathbf{a}^k=\oplus_{j\in A_k}a_j$.

 Since elements $f_j$ from $D(T_j)$ are continuous and thus
essentially bounded on $I_i^n$ for any $n\in\mathbb{N}$,
$\oplus_{j\in A_k}f_j \in \oplus_{j\in A_k}D(T_j)$ implies that
$$\mathrm{Range} \left[\oplus_{j\in
A_k}E^j\right](\Delta)\subseteq\oplus_{j\in A_k}D(T_j)\subset
\oplus_{j\in A_k}
    L^{\infty}(I_j^n)$$ and hence, we obtain
\begin{equation}\label{essbound}
\int_\Delta |W^k(\vee x_j, \lambda)|^2\,d\mu_{\mathbf{a}^k}\in
\oplus_{j\in A_k} L^{\infty}(I_j^{n})\,\,\, \forall n\in
\mathbb{N}.
\end{equation}
In \cite[XIII.5, Theorem 1(I)]{danford} it was shown that if we
have ordered representations $U_i$ of $L^2_i$ relative to the
operators $T_i$, $i\in\Omega$, the following formula is valid for
$j\in\Omega_1$:
$$
(U_jw_j)(\lambda)=\lim_{n\to\infty}\int_{I_j^n}\,w(x_j)\,\overline{W_j(x_j,
\lambda)}\,dx_j,
    \,\,\,\,\, w_j \in L^2_j,
$$
where the limit exists in $L^2(e_j,\mu_{a_j})$. Taking direct sums
in both sides of the last equality, for each system of compact
subintervals we obtain
$$
(U^k \oplus_{j\in A_k}w^n_j)(\lambda)=\oplus_{j\in
A_k}\int_{I_j^n}\,w_j(x_j)\,\overline{W_j(x_j, \lambda)}\,dx_j,
    \,\,\,\,\, w^n_j=w_j\chi_{I_j^n}.
$$
From (\ref{essbound}), it follows that for any bounded Borel set
$\Delta\in e_j$ and $I^{k,n}=\bigvee_{j\in A_k} I_j^n$,
$$\int_{I^{k,n}}\int_\Delta |W^k(\vee
x_j,\lambda)|^2\,d\mu_{\mathbf{a}^k}\,d(\vee x_j)<\infty$$ and
since $\mathbf{w}^k=\oplus_{j\in A_k}w_j(x_j)$ is assumed to
belong to $\oplus_{j\in A_k}L^2_j$, we may write:
$$
(U^k \mathbf{w}^{k,n})(\lambda)=\int_{I^{k,n}}\,\mathbf{w}^k(\vee
x_j)\,\overline{W^k(\vee x_j, \lambda)}\,d(\vee x_j).
$$
Taking the limit in the both sides and defining
$\mathbf{w}^{k}=\oplus_j{\lim_{n\to\infty}w_j^n}$ we obtain the
formula
\begin{equation}\label{untrans}
    (U^k \mathbf{w}^k)(\lambda)=\lim_{n\to\infty}\int_{I^{k,n}}\,\mathbf{w}^k(\vee x_j)\,\overline{W^k(\vee x_j, \lambda)}\,d(\vee
    x_j),
    \,\,\,\,\, \mathbf{w}^k\in \oplus_{j \in
    A_k}L^2_j.
\end{equation}
Note, that since for all $p=\overline{1,m_i}$ there exists the
equality $(\tau_i-\lambda)W_i^p=0$ (see \cite[XIII.5, Theorem
1]{danford}), it is obvious that $(\oplus_{j\in
A_k}\tau_j-\lambda)W^k=0$, where $\tau_j=\tau_i$ for a fixed $i$
and all $p=\overline{1,m_i}$. If $P(\epsilon(T_i))\cap
P(\epsilon(T_j))$ has zero spectral measures for all
$i,j\in\Omega$, then $A_k:\Omega_1=\cup_{k=1}^{\Lambda_1} A_k$ may
be constructed such that $A_k$ contains of indices
$\{i,k\},\,i\in\Omega,\,k=\overline{1,\max_i\{{m_i}\}}$.

\textbf{(B)} Consider the set of indices $\Omega_2=\{j\in\Omega_1:
j=\{i,1\},i\in\Omega\}$. Construct
$A_k:\Omega_2=\cup_{k=1}^{\Lambda_2} A_k$. Apply the reasonings
used in \textbf{(A)}, considering everywhere $\Omega_2$ instead of
$\Omega_1$. Hence, for each $A_k$ and we find a vector-function
$W_1^k(\vee x_j,\lambda)$ which is the solution of the equation
$(\oplus_{j\in A_k}\tau_j-\lambda)\mathbf{y}=0$. Consider $W_1^k$
and $W_1^s$ for $s\neq k$. For $\mathbf{a}^k$ there exists the
decomposition $\mathbf{a}^k = \mathbf{a}^k_k\oplus
\mathbf{a}^k_{k,s}$ (see the proof of Theorem \ref{ordered
representation}). This fact induces the decomposition for $W_1^k$:
$W_1^k=W^k_{1,k}\oplus W^k_{1,k,s}$. It is clear that being the
restrictions of $W_{1}^k$, the vector-functions $W^k_{1,k}$ and
$W^k_{1,k,s}$ are also the solutions of the equation
$(\oplus_{j\in A_k}\tau_j-\lambda)\mathbf{y}=0$. They, along with
$\mathbf{a}^k_k$ and $\mathbf{a}^k_{k,s}$ define unitary
transformations $U^k_k$ and $U_{k,s}^k$ by formula
(\ref{untrans}), such that:
$$U^k_k:\mathbf{L}^2(\mathbf{a}^k_{k})\to L^2(\mathbb{R},\mu_{k})\,\,\mbox{and}\,\,U_{k,s}^k:\mathbf{L}^2(\mathbf{a}^k_{k,s})\to
L^2(\mathbb{R},\mu_{k,s})$$ (see the definitions in the proof of
Theorem \ref{ordered representation}). This implies, that the
decomposition $W^k=W^k_{1,k}\oplus W^k_{1,k,s}$ is correct.

Define as $\max\{W^k_{1,k,s},W^s_{1,s,k}\}$ the vector-function,
which corresponds to the vector
$\max\{\mathbf{a}^k_{k,s},\mathbf{a}^s_{s,k}\}$, respectively
$\min\{W^k_{1,k,s}, W^s_{1,s,k}\}$ as the vector-function which
corresponds to that $\mathbf{a}^k_{k,s}$ or $\mathbf{a}^s_{s,k}$,
which is not maximal of the two.

\textbf{(C)} Without loss of generality, suppose that $k=1$ and
$s=2$. From the reasonings presented in Unit \textbf{(A)} of this
proof, it follows that
$$\Theta_1^{1\oplus 2}= W^1_{1,1}\oplus
W^2_{1,2}\oplus \max\left\{W^1_{1,1,2},W^2_{1,2,1}\right\}
$$ is correctly
constructed vector-function satisfying the statement of the
theorem for the case $T=[\oplus_{j\in A_1}T_j]\oplus [\oplus_{q\in
A_2}T_q]$. Apply the above described process to $\Theta_1^{1\oplus
2}$ and $W_1^3$ to obtain the correctly constructed
vector-function:
$$\Theta_1^{1\oplus 2\oplus 3}= \Theta^{1\oplus 2}_{1,1\oplus 2}\oplus
W^3_{1,3}\oplus \max\left\{\Theta^{1\oplus 2}_{1,1\oplus
2,3},W^3_{1,3,1\oplus 2}\right\}.
$$
Continuing this process, we finally obtain:
\begin{multline*}
\Theta_1(\vee x_i,\lambda)=\Theta_1^{1\oplus \cdots\oplus
\Lambda_2}=\\= \Theta^{1\oplus \cdots\oplus
\Lambda_2-1}_{1,1\oplus \cdots\oplus \Lambda_2-1}\oplus
W^{\Lambda_2}_{1,\Lambda_2}\oplus \max\left\{\Theta^{1\oplus
\cdots\oplus \Lambda_2-1}_{1,1\oplus \cdots\oplus
\Lambda_2-1,\Lambda_2},W^{\Lambda_2}_{1,\Lambda_2,1\oplus
\cdots\oplus \Lambda_2-1}\right\},
\end{multline*}
where in the case of $\Lambda_2=\infty$, $\Theta_1^{1\oplus
\cdots\oplus \Lambda_2}$ is the function which satisfies
(analogously to (\ref{convmaxv})):
\begin{multline}\label{converge1}
\left([\dir{E^i(\Delta)}]\,\int_\Delta\Theta_1^{1\oplus
\cdots\oplus
\Lambda_2}d\theta(\lambda),\int_\Delta\Theta_1^{1\oplus
\cdots\oplus
\Lambda_2}d\theta(\lambda)\right)=\\=\lim_{L\to\infty}
\left([\oplus_{j=1}^L{E^j(\Delta)}]\,\int_{\Delta}\Theta_1^{1\oplus
\cdots\oplus L}\,d\theta_L(\lambda),\int_{\Delta}\Theta_1^{1\oplus
\cdots\oplus L}\,d\theta_L(\lambda)\right),
\end{multline}
 for any
bounded Borel set $\Delta$, where
$\theta_L(\cdot)=\left([\oplus_{j=1}^L
E^j(\cdot)]\,\mathbf{a}^{1\oplus \cdots \oplus
L},\mathbf{a}^{1\oplus \cdots \oplus L}\right)$ is the measure of
the ordered representation of the space $\oplus_{j=1}^L L^2_j$.
The limit on the right side exists since for any bounded Borel
$\Delta$:
\begin{multline*}\left([\oplus_{j=1}^L{E^j(\Delta)}]\,\int_{\Delta}\Theta_1^{1\oplus
\cdots\oplus L}\,d\theta_L(\lambda),\int_{\Delta}\Theta_1^{1\oplus
\cdots\oplus L}\,d\theta_L(\lambda)\right)=
\\ =\left([\oplus_{j=1}^L{E^j(\Delta)}]\,\mathbf{a}^{1\oplus \cdots
\oplus L},\mathbf{a}^{1\oplus \cdots \oplus L}\right)\leqslant
\left([\oplus_{i=1}^{\infty}{E^i(\Delta)}]\,\oplus_{i=1}^{\infty}a_i,\oplus_{i=1}^{\infty}a_i\right)<\infty,
\end{multline*}
for all $L\in\mathbb{N}$ (Lemma \ref{irequality}). Despite seeming
weak, such convergence is quite natural. Indeed, (\ref{converge1})
implies that the cyclic subspace
$$\mathbf{L}^2\left(\int_{\Delta}\Theta_1^{1\oplus \cdots\oplus
L}\,d\theta_L(\lambda)\right)$$ is $\varepsilon$-close with the
cyclic subspace $$\mathbf{L}^2\left(\int_{\Delta}\Theta_1^{1\oplus
\cdots\oplus \Lambda_2}\,d\theta(\lambda)\right),$$ when $L$ is
sufficiently big. That is, in the topology of $ \mathbf{L}^2$ for
any Borel set $\Delta$,
$$f(T)\int_{\Delta}\Theta_1^{1\oplus \cdots\oplus
L}\,d\theta_L(\lambda)\to f(T)\int_{\Delta}\Theta_1^{1\oplus
\cdots\oplus \Lambda_2}\,d\theta(\lambda),$$ for any Borel $f$ as
$L\to\infty$. This means that
$$\int_{\Delta}\Theta_1^{1\oplus \cdots\oplus
L}\,d\theta_L(\lambda)\to \int_{\Delta}\Theta_1^{1\oplus
\cdots\oplus \Lambda_2}\,d\theta(\lambda),\,\mbox{as}\,\, L\to
\infty.
$$

 \textbf{(D)} Define
$\Omega_3=\{j\in\Omega_1:j=\{i,2\},i\in\Omega\}$. Construct
$A_k:\Omega_3=\cup_{k=1}^{\Lambda_3} A_k$. Apply processes
\textbf{(B)} and \textbf{(C)} of this proof, substituting
everywhere $\Omega_3$ instead of $\Omega_2$. We obtain a
vector-function $ \Theta_2^{1\oplus \cdots\oplus \Lambda_3}$,
which is defined on the set $\cup_i P(e_2^i)$. But, as we know
(see (\ref{mult-set})), the set $s_2$ also includes the sets where
there are non-empty superpositions of $\epsilon(T_i)$. Therefore,
designating
\begin{multline*}
\Theta_2^1=\Theta_2^{1\oplus \cdots\oplus \Lambda_3},\,
\Theta_2^2=\min\{W^1_{1,1,2},W^2_{1,2,1}\},\dots,\\
\Theta_2^{\Lambda_2+1}=\min\left\{\Theta^{1\oplus \cdots\oplus
\Lambda_2-1}_{1,1\oplus \cdots\oplus
\Lambda_2-1,\Lambda_2},W^{\Lambda_2}_{1,\Lambda_2,1\oplus
\cdots\oplus \Lambda_2-1}\right\},
\end{multline*}
we may again use the process \textbf{(C)} to build the
vector-function $\Theta_2(\vee x_i,\lambda)$ defined on $s_2$ and
$\Theta_2(\vee x_i,\lambda)=0$ for $\lambda\in \mathbb{R}\setminus
s_2$. Using processes \textbf{(B)}, \textbf{(C)}, \textbf{(D)} and
formula (\ref{mult-set}), we finally obtain $\Theta_m(\vee
x_i,\lambda)$.

\textbf{(E)} The above presented constructions show, that all
vector-functions $\Theta_k(\vee x_i,\lambda)$, $k=\overline{1,m}$
are the solutions of the equation $(\oplus_{i\in
\Omega}\tau_i-\lambda)\mathbf{y}=0$, moreover they equal zero on
$\mathbb{R}\setminus s_k$ and satisfy formulas (\ref{essbound1})
and (\ref{untrans1}).

The last thing is to prove the linear independence. In order to
make the reasonings more transparent, we prove the linear
independence for the special case of two vector-functions
$$\Theta_1 = W^1_{1,1}\oplus W^2_{1,2}\oplus
\max\left\{W^1_{1,1,2},W^2_{1,2,1}\right\} $$ and
$$\Theta_2=\min\left\{W^1_{1,1,2},W^2_{1,2,1}\right\}.$$ Without
loss of generality suppose that
$\max\left\{W^1_{1,1,2},W^2_{1,2,1}\right\}=W^1_{1,1,2}$. It is
clear that
\begin{multline*}
\alpha\Theta_1+\beta\Theta_2=\\=\alpha
\left(\{W^1_{1,1},0,0,0\}+\{0,0,W^2_{1,2},0\}+\{0,W^1_{1,1,2},0,0\}\right)+
\beta \{0,0,0,W^2_{1,2,1}\}=\\=\{\alpha W^1_{1,1},\alpha
W^1_{1,1,2},\alpha W^2_{1,2}, \beta W^2_{1,2,1}\} =0
\end{multline*}
implies $\alpha=\beta=0$. The linear independence in the general
case is proved using the same ideas. Thus, the linear independence
is proved and this finishes the proof of the theorem.
\end{proof}
Note, that the given proof introduces the general method of
constructing eigenfunctions for a vector-operator. For theoretical
purposes, the form of the obtained eigenfunctions could be
simplified by totally ordering the set
$\{T^j\}_{j=1}^{\Lambda_2}$. This is achieved by saying that
$T^k\preceq T^s$ if $\max\{W^k_{1,k,s},W^s_{1,s,k}\}=W^k_{1,k,s}$.
At that, $T^k\simeq T^s$ if and only if $T^k\preceq T^s$ and
$T^s\preceq T^k$. According to this, we build
$\oplus_{j=1}^{\Lambda_2} T^j$, where $T^j\preceq T^{j+1}$,
$j=\overline{1,\Lambda_2-1}$ if $\Lambda_2\geqslant 2$. The
obtained vector-operator is obviously equivalent to the initial
vector-operator (comprising unordered operators).

As an important corollary of Theorem \ref{eigexpthm} we obtain
\begin{thm}[Eigenfunction expansions]
For any $\mathbf{w}\in \mathbf{L}^2$, there exists a decomposition
\begin{equation*}\label{eigenfunction expansion}
\mathbf{w}= \sum_{k=1}^m\,\lim_{n\to\infty} \int_{-n}^{+n}(U
\mathbf{w})^k(\lambda)\Theta_k(\vee x_i,
\lambda)\,d\theta(\lambda).
\end{equation*}
\end{thm}

\begin{proof} From the process of building $\Theta_k$ in the previous
proof and, in particular (\ref{ran}), it follows that
\begin{equation*}\label{ran2}
    \left(E(\Delta)F(T)\mathbf{a}^k\right)(\vee x_i)=\int_\Delta
    \Theta_k(\vee
    x_i,\lambda)F(\lambda)\,d\mu_{\mathbf{a}_k}(\lambda).
\end{equation*}
Substituting here $F=(U\mathbf{w})^k$, we obtain
\begin{multline*}
\int_{-n}^n
    \Theta_k(\vee
    x_i,\lambda)(U\mathbf{w})^k(\lambda)\,d\mu_{\mathbf{a}_k}(\lambda)= E[-n,n]F(T)\mathbf{a}^k
   \to F(T)\mathbf{a}^k = {U^{k}}^{-1}F = \mathbf{w}^k.
\end{multline*}
Now the statement of the theorem becomes clear, since
$\mathbf{w}=\oplus_{k=1}^m{\mathbf{w}_k}$.
\end{proof}

\begin{thm}\label{intoper}
Let $T$ be a self-adjoint vector-operator, generated by an EMZ
system $\{I_i,\tau_i\}_{i\in\Omega}$. Let the measure $\theta$ and
the sets $\{s_k\}_{k=1}^m$ be respectively a measure and
multiplicity sets of an ordered representation of the space
$\mathbf{L}^2=\oplus_{i\in\Omega}L^2(I_i)$, relative to the
operator $T$. The kernels $\{\Theta_k\}_{k=1}^m$ are the
generalized vector-operator eigenfunctions, corresponding to the
multiplicity sets (as defined in Theorem \ref{eigexpthm}). Given a
bounded Borel function $F$, which equals zero beyond a compact
Borel set $\Delta$, the bounded vector-operator $F(T)$ may be
represented as an integral operator:
\begin{equation}\label{intoperator}[F(T)\mathbf{f}](\vee
s_i)=\int_{I} \mathbf{f}(\vee x_i)K(F;\vee s_i, \vee x_i)\, d(\vee
x_i),\end{equation} where $\mathbf{f}\in \mathbf{L}^2$ and
\begin{equation}\label{kernel} K(F;\vee x_i, \vee s_i) = \sum_{k=1}^m \int_\Delta
F(\lambda) \Theta_k(\vee x_i,\lambda)\overline{\Theta_k(\vee
s_i,\lambda)}\,d\theta(\lambda).\end{equation}
\end{thm}
\begin{proof}
From Theorem \ref{eigenfunction expansion} it follows that
$$[F(T)\mathbf{f}](\vee
x_i) = \sum_{k=1}^m\,\int_\Delta(U
F(T)\mathbf{f})^k(\lambda)\Theta_k(\vee x_i,
\lambda)\,d\theta(\lambda),$$ and since for any spectral
representation $$(U
F(T)\mathbf{f})^k(\lambda)=F(\lambda)(U\mathbf{f})^k(\lambda),$$
from Theorem \ref{eigexpthm} we obtain:
\begin{multline}\label{fott}
[F(T)\mathbf{f}](\vee s_i) = \sum_{k=1}^m\,\int_\Delta
F(\lambda)(U\mathbf{f})^k(\lambda)\Theta_k(\vee x_i,
\lambda)\,d\theta(\lambda)=\\= \sum_{k=1}^m\,\int_\Delta
F(\lambda)\Theta_k(\vee s_i, \lambda)\int_{I}\,\mathbf{f}(\vee
x_i)\,\overline{\Theta_k(\vee x_i, \lambda)}\,d(\vee
    x_i)\,d\theta(\lambda),
\end{multline}
where $$\int_{I}\,\mathbf{f}(\vee x_i)\,\overline{\Theta_k(\vee
x_i, \lambda)}\,d(\vee
    x_i)=\lim_{n\to\infty}\int_{I^{n}}\,\mathbf{f}(\vee x_i)\,\overline{\Theta_k(\vee x_i, \lambda)}\,d(\vee
    x_i),$$ for which see formula (\ref{untrans1}).

Note that (see \cite[XII.3.8]{danford})
$F(T)\mathbf{f}\in\oplus_{i\in\Omega}\left(\cap_{n=1}^\infty
D(T_i^n)\right).$ For any system $\{J_i\}_{i\in\Omega}$ of compact
subintervals of the respective intervals from
$\{I_i\}_{i\in\Omega}$, define the space
$\oplus_{i\in\Omega}C(J_i)=C(J)$, $J=\bigvee_iJ_i,$ as the space
of continuous vector-functions with the norm
$$\|\mathbf{f}\|_{C(J)}=\sup_i \sup_{s_i\in
J_i}|f_i(s_i)|.$$ Hence, the mapping $\mathbf{f}\to
F(T)\mathbf{f}$ is continuous as the operator from $\mathbf{L}^2$
to $C(J)$. This means that there exists a constant $M(J)$, such
that
$$\|F(T)\mathbf{f}\|_{C(J)}\leqslant
M(J)\|\mathbf{f}\|_{\mathbf{L}^2},$$ or
\begin{equation}\label{supremumc}\sup_i
\sup_{s_i\in J_i}|(F(T_i)f_i)(s_i)|\leqslant
M(J)\|\mathbf{f}\|_{\mathbf{L}^2}.
\end{equation}
Let $m<\infty$. For each $i\in \Omega$, define  $\mathcal{H}_i$ as
a dense set in $L^2_i$ consisting of functions equalling zero
beyond a compact subset of $I_i$. We can interchange the integrals
in (\ref{fott}):
\begin{equation}\label{densesusetformula}[F(T)\mathbf{f}](\vee s_i)=\int_{I}
\mathbf{f}(\vee x_i)K(F,\vee x_i, \vee s_i)\, d(\vee
x_i),\end{equation} for $\mathbf{f}\in \oplus_{i\in\Omega}
\mathcal{H}_i$ and $$ K(F;\vee x_i, \vee s_i) = \sum_{k=1}^m
\int_\Delta F(\lambda) \Theta_k(\vee
x_i,\lambda)\overline{\Theta_k(\vee
s_i,\lambda)}\,d\theta(\lambda). $$

From (\ref{supremumc}) we obtain:
$$\sup_i\sup_{s_i\in J_i}\left|\int_{I}
\mathbf{f}(\vee x_i)K(F,\vee x_i, \vee s_i)\, d(\vee
x_i)\right|\leqslant M(J)\|\mathbf{f}\|_{\mathbf{L}^2}.$$  It is
clear now, that the formula (\ref{densesusetformula}) holds for
any $\mathbf{f}\in\mathbf{L}^2$.

Pass now to the case $m=\infty$. Recall that $$U:\mathbf{L}^2\to
\oplus_{k=1}^\infty L^2(s_k,\theta).$$ For each $n<\infty$ define
an orthogonal projector $P_n:\mathbf{L}^2\to\mathbf{L}^2$, such
that $$P_n \oplus_{i\in\Omega}f_i = \{f_1,\dots,f_n,0,0,\dots\}.$$
Define  continuous linear functionals
$\phi_{s_i}(f_i)=(F(T_i)f_i)(s_i)$, for which there exist
$g_{s_i}\in L^2_i$ such that $(F(T_i)f_i)(s_i)=(f_i,g_{s_i})_i$.
From the reasonings presented in the beginning of this proof for
the case of the finite multiplicity, we obtain:
\begin{multline*}(F(T)U^{-1}P_nU\mathbf{f})(\vee s_i)=(U^{-1}P_nUF(T)\mathbf{f})(\vee
s_i)=\\=\sum_{k=1}^n\,\int_\Delta F(\lambda)\Theta_k(\vee s_i,
\lambda)\int_{I}\,\mathbf{f}(\vee x_i)\,\overline{\Theta_k(\vee
x_i, \lambda)}\,d(\vee
    x_i)\,d\theta(\lambda).\end{multline*} That is  \begin{equation}\label{resolutionfinite}(F(T)U^{-1}P_nU\mathbf{f})(\vee s_i)=\int_{I}
\mathbf{f}(\vee x_i)K_n(F;\vee x_i, \vee s_i)\, d(\vee
x_i),\end{equation} where
$$K_n(F;\vee x_i, \vee s_i) = \sum_{k=1}^n \int_\Delta F(\lambda)
\Theta_k(\vee x_i,\lambda)\overline{\Theta_k(\vee
s_i,\lambda)}\,d\theta(\lambda).$$ Since
\begin{multline*}(F(T)U^{-1}P_nU\mathbf{f})(\vee
s_i)=\oplus_{i\in\Omega}(F(T_i)[U^{-1}P_nU]_if_i)(s_i)=\\=
\left\{([U^{-1}P_nU]_1f_1,g_{s_1})_1,
([U^{-1}P_nU]_2f_2,g_{s_2})_2,\dots,
([U^{-1}P_nU]_jf_j,g_{s_j})_j,\dots\right\}=\\=
\left\{(f_1,[U^{-1}P_nU]_1g_{s_1})_1,
(f_2,[U^{-1}P_nU]_2g_{s_2})_2,\dots,
(f_j,[U^{-1}P_nU]_jg_{s_j})_j,\dots\right\},\end{multline*}since
the coordinate of a unitary vector-operator is unitary in the
coordinate space. From this formula and (\ref{resolutionfinite})
we obtain that the coordinate $K^i_n(F;x_i, \cdot)$ of $K_n(F;\vee
x_i, \cdot)$ satisfies the equation $$\overline{K^i_n(F;x_i,
s_i)}=[U^{-1}P_nU]_ig_{s_i}.$$ From  which follows that
$$\lim_{n\to \infty}\overline{K^i_n(F;x_i,
s_i)}=\lim_{n\to\infty}[U^{-1}P_nU]_ig_{s_i},$$ and thus
$\overline{K^i(F;x_i, s_i)}=g_{s_i}$ and $\overline{K(F;\vee x_i,
\vee s_i)}=\oplus_{i\in\Omega}g_{s_i}$ which means that the series
defining $K(F;x_i, s_i)$ converges in $\mathbf{L^2}$ for each
fixed $\vee s_i$. Moreover, \begin{multline*}\int_{I}
\mathbf{f}(\vee x_i)K(F;\vee x_i, \vee s_i)\, d(\vee
x_i)=\oplus_{i\in\Omega}\int_{I} f_i(x_i)K^i(F; x_i,  s_i)\, dx_i
= \\ = \left\{(f_1,g_{s_1})_1, (f_2,g_{s_2})_2,\dots,
(f_j,g_{s_j})_j,\dots\right\} = \\=\oplus_{i\in\Omega}
(F(T_i)f_i)(s_i) = (F(T)\mathbf{f})(\vee s_i).
\end{multline*}

The theorem is proved.
\end{proof}

Since the kernels from Theorem \ref{eigexpthm} are only measurable
relative to $\lambda$, the following theorem is important to
strengthen the practical value of Theorems \ref{eigexpthm},
\ref{eigenfunction expansion} and \ref{intoper}:

\begin{thm}\label{kernelsdecomp}Each kernel $\Theta_k(\vee
x_i,\lambda)$, $k=\overline{1,m}$, may be decomposed as
\begin{equation}\label{analdec}
\Theta_k(\vee x_i,\lambda)=\sum_{s=1}^{M_k}
\gamma_{sk}(\lambda)\sigma_{sk}(\vee x_i,\lambda),
\end{equation}
where the $M_k$ are finite for each $k$ and $\sigma_{sk}(\vee
x_i,\lambda)$ depend analytically on $\lambda$.

\end{thm}
\begin{proof}
We separate the proof in parts which will correspond to the
analogous parts of the proof of Theorem  \ref{eigexpthm}.

\textbf{(A')} Each kernel $W_j(x_j,\lambda)$ from the part
\textbf{(A)} of the proof of Theorem \ref{eigexpthm} may be
decomposed:
$$W_j(\cdot,\lambda)=\sum_{s=1}^{n_j}\alpha_{js}(\lambda)\sigma_{js}(\cdot,\lambda),$$
where $\alpha_{js}$ are supposed to equal zero on
$\mathbb{R}\setminus e_j$, see \cite[p. 1351]{danford}.
Supplementing the defining systems with zeros where necessary, we
obtain:
\begin{multline*}
W^k(\vee x_j,\lambda)= \oplus_{j\in A_k} W_j(x_j,\lambda) =
\oplus_{j\in A_k}
\sum_{s=1}^{n_j}\alpha_{js}(\lambda)\sigma_{js}(x_j,\lambda)=\\=
\sum_{q=1}^{N_k}\oplus_{j\in
A_k}\alpha_{jq}(\lambda)\sigma_{jq}(x_j,\lambda)=\sum_{q=1}^{N_k}\alpha^k_{q}(\lambda)\sigma^k_{q}(\vee
x_j,\lambda),
\end{multline*}
where $$N_k=\max_{j\in A_k} n_j,\,\,
\alpha^k_{q}(\lambda)=\sum_{j\in A_k}\alpha_{jq}(\lambda)\,\,
\mbox{and}\,\, \sigma^k_{q}(\vee x_j,\lambda)=\oplus_{j\in
A_k}\sigma_{jq}(x_j,\lambda).$$

Since $e_j$ and $e_k$ do not intersect almost everywhere for
$j,k\in\Omega_2$, $j\neq k$, the series $\sum_{j\in
A_k}\alpha_{jq}(\lambda)$ converges almost everywhere on
$\cup_{j\in\Omega_2} P(e_j)$.

\textbf{(C')} Now pass to the part \textbf{(B)}. There we obtained
the decompositions $W_1^k=W^k_{1,k}\oplus W^k_{1,k,s}$ and
$W_1^s=W^s_{1,s}\oplus W^s_{1,s,k}$. Let us totally order the set
$\{T^j\}_{j=1}^{\Lambda_2}$ saying that $T^k\preceq T^s$ if
$\max\{W^k_{1,k,s},W^s_{1,s,k}\}=W^k_{1,k,s}$. At that, $T^k\simeq
T^s$ if and only if $T^k\preceq T^s$ and $T^s\preceq T^k$.
According to this, we build $\oplus_{j=1}^{\Lambda_2} T^j$, where
$T^j\preceq T^{j+1}$, $j=\overline{1,\Lambda_2-1}$ if
$\Lambda_2\geqslant 2$. The obtained vector-operator is obviously
equivalent to the initial vector-operator (comprising unordered
operators). Note that $$W_1^k(\vee
x_i,\lambda)=\sum_{q=1}^{N_k}\alpha^k_{1q}(\lambda)\sigma^k_{1q}(\vee
x_j,\lambda)$$ and analogously $$W_1^s(\vee
x_i,\lambda)=\sum_{p=1}^{N_s}\alpha^s_{1p}(\lambda)\sigma^s_{1p}(\vee
x_j,\lambda).$$

All the above leads to the following:
\begin{multline*}
\Theta_1^{1\oplus 2}= W^1_{1,1}\oplus W^2_{1,2}\oplus
\max\left\{W^1_{1,1,2},W^2_{1,2,1}\right\} = W_1^1\oplus W^2_{1,2}
= \\ =
\left(\sum_{q=1}^{N_1}\alpha^1_{1q}(\lambda)\sigma^1_{1q}(\vee
x_j,\lambda)\right)\oplus
\left(\sum_{p=1}^{N_2}\alpha^2_{1p}(\lambda)\chi_{\epsilon_{2}}(\lambda)\sigma^2_{1p}(\vee
x_j,\lambda)\right) =\\= \sum_{s=1}^{N^{1\oplus 2}}\alpha^{1\oplus
2}_{1s}(\lambda)\sigma^{1\oplus 2}_{1s}(\vee x_j,\lambda)
\end{multline*}
where $$N^{1\oplus 2}=\max\{N_1,N_2\};\,\,\alpha^{1\oplus
2}_{1s}(\lambda)=\alpha^1_{1s}(\lambda)+\alpha^2_{1s}(\lambda)\chi_{\epsilon_{2}}(\lambda),$$$$\sigma^{1\oplus
2}_{1s}(\vee x_j,\lambda)=\sigma^1_{1s}(\vee
x_j,\lambda)\oplus\left(\sigma^2_{1s}(\vee
x_j,\lambda)\chi_{\epsilon_{2}}(\lambda)\right),\,\,s=\overline{1,N^{1\oplus
2}}.
$$

Continuing this process till the finite $\Lambda_2$, we obtain:
\begin{equation}\label{def1}
 \Theta_1(\vee x_i,\lambda)=\Theta_1^{1\oplus \cdots\oplus
\Lambda_2}=
\sum_{s=1}^{N^{1\oplus\cdots\oplus\Lambda_2}}\alpha^{1\oplus\cdots\oplus\Lambda_2}_{1s}(\lambda)\sigma_{1s}^{1\oplus\cdots\oplus\Lambda_2}(\vee
x_j,\lambda),
\end{equation}
where $N^{1\oplus \cdots\oplus \Lambda_2}=\max\{N_1,N_2,\dots
N_{\Lambda_2}\}$ and for $s=\overline{1,N^{1\oplus \cdots\oplus
\Lambda_2}}$:
\begin{equation}\label{coef1}
\alpha^{1\oplus \cdots\oplus
\Lambda_2}_{1s}(\lambda)=\alpha^1_{1s}(\lambda)+\sum_{i=2}^{\Lambda_2}\alpha^i_{1s}(\lambda)\chi_{\epsilon_{i}}(\lambda);
\end{equation}

$$\sigma^{1\oplus
\cdots\oplus \Lambda_2}_{1s}(\vee x_j,\lambda)=\sigma^1_{1s}(\vee
x_j,\lambda)\oplus\left(\oplus_{i=2}^{\Lambda_2}\sigma^i_{1s}(\vee
x_j,\lambda)\chi_{\epsilon_{i}}(\lambda)\right).
$$
In the case of infinite $\Lambda_2$, $N^{1\oplus \cdots\oplus
\Lambda_2}$ is clearly finite. The series in the right side of
(\ref{coef1}) pointwise converges, since it consists of items
defined on non-intersecting sets. $\sigma^{1\oplus \cdots\oplus
\Lambda_2}_{1s}(\vee x_j,\lambda)$ is defined by induction as the
element which satisfies (see \ref{converge1})
$$\lim_{L\to\infty}\left\|\int_\Delta\sigma^{1\oplus \cdots\oplus
\Lambda_2}_{1s}(\vee
x_j,\lambda)\,d\theta(\lambda)-\int_\Delta\sigma^{1\oplus
\cdots\oplus L}_{1s}(\vee
x_j,\lambda)\,d\theta_L(\lambda)\right\|=0.$$ Since for each
finite iteration the equality (\ref{def1}) is fulfilled, it is
clear that for an infinite $\Lambda_2$ it will be fulfilled too.

\textbf{(D')} Borrowing the designations from \textbf{(D)} and
using processes described in \textbf{(A')} and \textbf{(C')}, we
shall come to the decomposition of $\Theta_2^{1\oplus \cdots\oplus
\Lambda_3}$:

\begin{equation*}
\Theta_2^{1\oplus \cdots\oplus \Lambda_3}=
\sum_{s=1}^{N^{1\oplus\cdots\oplus\Lambda_3}}\alpha^{1\oplus\cdots\oplus\Lambda_3}_{2s}(\lambda)\sigma_{2s}^{1\oplus\cdots\oplus\Lambda_3}(\vee
x_j,\lambda).
\end{equation*}

To obtain $\Theta_2(\vee x_i,\lambda)$, as in \textbf{(D)}, we
repeat part \textbf{(C')} for
$$
\Theta_2^1=\Theta_2^{1\oplus \cdots\oplus \Lambda_3},\,
\Theta_2^2=W^2_{1,2,1},\dots,
\Theta_2^{\Lambda_2+1}=W^{\Lambda_2}_{1,\Lambda_2,1\oplus
\cdots\oplus \Lambda_2-1}.
$$
Finally, the same way we obtain decompositions for all
$\Theta_k(\vee x_i,\lambda)$, $k=\overline{1,m}$, which will have
the form (\ref{analdec}).
\end{proof}

In Theorem \ref{kernelsdecomp} we have obtained the formula
(\ref{analdec}) which we substitute now in (\ref{kernel}). Thus we
obtain:
$$K(F;\vee x_i, \vee s_i) = \sum_{k=1}^m \int_\Delta
 \sum_{s,p=1}^{M_k}F(\lambda) \gamma_{sk}(\lambda)\overline{\gamma_{pk}(\lambda)}\sigma_{sk}(\vee
x_i,\lambda)\overline{\sigma_{pk}(\vee
s_i,\lambda)}\,d\theta(\lambda).$$ The last formula may be
rewritten as
$$K(F;\vee x_i, \vee s_i) = \int_\Delta
 \left\{\sum_{s,p=1}^{M_k}F(\lambda) \sigma_{sk}(\vee
x_i,\lambda)\overline{\sigma_{pk}(\vee
s_i,\lambda)}\,d\varrho_{sp}(\lambda)\right\},$$ where
$\varrho_{sp}(\lambda)=\sum_{k=1}^m
\int_\Delta\gamma_{sk}(\lambda)\overline{\gamma_{pk}(\lambda)}\,d\theta(\lambda).$
Let us verify that
$\mathcal{R}_{M_k}(\lambda)=\{\varrho_{sp}(\lambda)\}$ is a
correctly constructed matrix measure. First of all if $\{\xi_1,
\dots \xi_{M_k}\}$ is a collection of complex numbers, then
$$\sum_{s,p=1}^{M_k}\varrho_{sp}(\Delta)\xi_s\overline{\xi_p}=\sum_{k=1}^m \int_\Delta \left|\sum_{s=1}^{M_k}\gamma_{ks}\xi_s\right|^2d\theta(\lambda)\geqslant0.$$
Moreover it is obvious that
$\mathcal{R}_{M_k}(\cup_i\Delta_i)=\sum_i\mathcal{R}_{M_k}(\Delta_i),$
where $\cup_i\Delta_i$ is precompact.

\textbf{Acknowledgements}. 1. The work has been carried out with
the financial support of the Abdus Salam International Center for
Theoretical Physics (grant AC-84). 2. The author is grateful to
Professor R.R.~Ashurov for his productive pieces of advice and
deep attention to the work.

\end{document}